\documentclass[conference] {IEEEtran}
\usepackage{amsthm}
\usepackage{amsmath}
\usepackage{amsfonts}
\usepackage{graphicx}
\newtheorem{thm}{Theorem}

\newtheorem{lemma}[thm]{Lemma}
\newtheorem{cor}[thm]{Corollary}
\theoremstyle{definition}

\newtheorem*{prf}{Proof}
\usepackage{theoremref}
\newcommand{\C}[0]{\mathbb{C}}
\newcommand{\R}[0]{\mathbb{R}}
\title{A Delsarte-Style Proof of the Bukh--Cox Bound}

\author{Mark Magsino, Dustin G. Mixon, Hans Parshall
        \\ Department of Mathematics
        \\ The Ohio State University} 
\begin{document}
\maketitle

\begin{abstract}
The line packing problem is concerned with the optimal packing of points in real or complex projective space so that the minimum distance between points is maximized. Until recently, all bounds on optimal line packings were known to be derivable from Delsarte's linear program. Last year, Bukh and Cox introduced a new bound for the line packing problem using completely different techniques. In this paper, we use ideas from the Bukh--Cox proof to find a new proof of the Welch bound, and then we use ideas from Delsarte's linear program to find a new proof of the Bukh--Cox bound. Hopefully, these unifying principles will lead to further refinements.
\end{abstract}
\section{Introduction}
The last decade has seen a surge of progress in \textit{the line packing 
problem}, where the objective is to pack $n$ points in 
$\mathbb{R}\mathbf{P}^{d-1}$ or $\mathbb{C}\mathbf{P}^{d-1}$ so that the 
minimum distance is maximized.
Indeed, while instances of this problem date back to Tammes \cite{Tam1930} 
and Fejes T\'{o}th \cite{Tot1965}, the substantial progress in recent 
days has been motivated by emerging applications in compressed sensing 
\cite{BanFicMixWon2013}, digital fingerprinting \cite{MixQuiKiyFic2013}, 
quantum state tomography 
\cite{RenBluScoCav2004}, and multiple description coding \cite{StrHea2003}.
Most progress in this direction has come from identifying new packings that 
achieve equality in the so-called Welch bound (see \cite{FicMix2015} for a 
survey), but last year, Bukh and Cox \cite{BukCox2018}
discovered a completely different bound, along with a large family of 
packings that achieve equality in their bound.

Focusing on the complex case, let $X=\{x_i\}_{i\in[n]}$ be a sequence of unit 
vectors in $\mathbb{C}^d$, and define coherence to be
\[
\mu(X)
:=\max_{1\leq i<j\leq n}|\langle x_i,x_j\rangle|.
\]
Then the Buhk--Cox bound guarantees
\[
\mu(X)
\geq \frac{(n-d)^2}{n+(n^2-nd-n)\sqrt{1+n-d}-(n-d)^2}
\]
provided $n>d$.
The Bukh--Cox bound is the best known lower bound on coherence in the 
regime where $n=d+O(\sqrt{d})$.
While the other lower bounds can be proven using Delsarte's linear
program \cite{DelGoeSei1991}, the proof of the Bukh--Cox bound is completely 
different: it hinges on an upper bound on the first moment of 
isotropic measures.

In the present paper, we provide an alternate proof of the Bukh--Cox bound.
We start by isolating a lemma of Bukh and Cox that identifies a fundamental 
duality between the coherence of $n=d+k$ unit vectors in $d$ dimensions and 
the entrywise $1$-norm of the Gram matrix of $\frac{n}{k}$-tight frames of 
$n$ vectors in $k$ dimensions.
Next, we illustrate the power of this lemma by using it to find a new proof 
of the Welch bound.
Finally, we combine the lemma with ideas from Delsarte's linear program to 
obtain a new proof of the Bukh--Cox bound.
This new proof helps to unify the existing theory of line packing, and 
hopefully, it will spur further improvements (say, by leveraging ideas 
from semidefinite programming \cite{LaaVal2015}).

\section{The Bukh--Cox Lemma}
Let $X = \{x_i\}_{i=1}^n$ denote any sequence in $\mathbb{C}^d$.
By abuse of notation, we identify $X$ with the $d\times n$ matrix whose $i$th column is $x_i$.
We say $X$ is a \textit{$c$-tight frame} if $XX^*=cI$.
Let $N(d,n)$ denote the set of matrices in $\mathbb{C}^{d\times n}$ with unit norm columns, and let $T(d,n)$ denote the set of matrices in $\mathbb{C}^{d\times n}$ corresponding to $\frac{n}{d}$-tight frames.
Define
\[
\gamma(d,n) := \max_{Y \in T(d,n)} \|Y^* Y\|_1.
\]
(Indeed, the maximum exists by a compactness argument.)
We say $X\in N(d,n)$ is \textit{equiangular} if there exists a constant $c$ such that
$|\langle x_i, x_j\rangle| = c$ for every $1\leq i<j\leq n$. 
With this nomenclature, we are ready to state the following lemma of Bukh and Cox:

\begin{lemma}\thlabel{thm: minmax coherence}
Let $n = d+k$. Then every $X \in N(d,n)$ satisfies
\begin{equation}\label{eqn: minmax coherence}
    \mu(X) \geq \frac{n}{\gamma(k,n) - n}.
\end{equation}
Furthermore, $X$ minimizes $\mu(X)$ over $N(d,n)$ if $X$ is equiangular
and there exists $Y =\{y_i\}_{i=1}^n\in T(k,n)$ such that 
\begin{enumerate}
    \item[(i)] $Y$ maximizes $\|Y^* Y\|_1$ over $T(k,n)$,
    \item[(ii)] $XY^* = 0$, and
    \item[(iii)] $\mathrm{sgn} \langle x_i, x_j \rangle = - \mathrm{sgn} 
          \langle y_i, y_j \rangle$ for
          $1 \leq i < j \leq n$.
\end{enumerate}
\end{lemma}

\begin{prf}
Given $X \in N(d,n)$, select $Y \in T(k,n)$ such that $XY^* = 0$. Then,
\begin{align}
0 &= (X^*X Y^* Y)_{ii} = \sum_{j=1}^n \langle x_i, x_j \rangle
    \langle y_j, y_i \rangle \nonumber \\
    ~ & = \|y_i\|_2^2 + \sum_{\substack{j=1 \\ j \neq i}}^n \langle x_i, x_j
    \rangle \langle y_j, y_i \rangle \label{eqn: double inner product}.
\end{align}
Bringing $\|y_i\|_2^2$ to the left hand side of 
(\ref{eqn: double inner product}) and taking absolute values, we have that
\begin{align}
    \|y_i\|_2^2 & = \bigg| \sum_{\substack{j=1 \\ j \neq i}}^n
    \langle x_i, x_j \rangle \langle y_j, y_i \rangle \bigg| 
    \leq \sum_{\substack{j=1 \\ j \neq i}}^n
    | \langle x_i, x_j \rangle| |\langle y_i, y_j \rangle| 
    \label{eqn: triangle ineq}\\
    ~& \leq  \mu(X)
    \sum_{\substack{j=1 \\ j \neq i}}^n | \langle y_i, y_j \rangle |,
    \label{eqn: coherence ineq} 
    \end{align}
where (\ref{eqn: triangle ineq}) uses the triangle inequality and (\ref{eqn:
coherence ineq}) is by the definition of coherence.
Using (\ref{eqn: coherence ineq}) and $YY^* = (n/k) I$, 
\begin{align}
    n &= \mathrm{tr}(YY^*) = \mathrm{tr}(Y^* Y) = \sum_{i=1}^n \|y_i\|_2^2 
    \nonumber \\ 
    ~& \leq \sum_{i=1}^n  \mu(X)
    \sum_{\substack{j=1 \\ j \neq i}}^n | \langle y_i,
    y_j \rangle | = \mu(X) (\|Y^* Y\|_1 - \mathrm{tr}(Y^* Y)) \nonumber \\
    ~& = \mu(X) (\|Y^* Y\|_1 - n) \nonumber.
\end{align}
Thus, we conclude that
\begin{equation}\label{eqn: max inequality}
    \mu(X) \geq \frac{n}{\|Y^* Y\|_1 - n} \geq \frac{n}{\gamma(k,n)
    - n}.
\end{equation}
This proves the bound.
Considering \eqref{eqn: double inner product}, equality occurs in \eqref{eqn: triangle ineq}, \eqref{eqn: coherence ineq} and \eqref{eqn: max inequality} when $X$ is equiangular and (i)--(iii) holds. \qed
%
\end{prf}

\section{The Welch Bound}
\begin{thm}\thlabel{lemma: welch y1}
For all $Y \in T(k,n)$ we have 
\begin{equation}\label{eqn: welch y}
    \|Y^* Y\|_1 \leq n + \left[ n(n-1)\left(\frac{n^2}{k} - n \right) 
    \right]^{1/2}.
\end{equation}
Equality is achieved if and only if $Y$ is an equiangular tight frame.
\end{thm}

\begin{prf}
First we separate the diagonal part of $\|Y^* Y\|_1$ and use the Cauchy--Schwarz
inequality,
\begin{align}
    \|Y^* Y\|_1 &= n + \sum_{i=1}^n \sum_{\substack{j=1 \\ j \neq i}}^n 
    |\langle y_i, y_j \rangle | \nonumber \\
    ~ & \leq n + \bigg[n(n-1) \sum_{i=1}^n \sum_{\substack{j=1 \\ j \neq i}}^n
    |\langle y_i, y_j\rangle|^2
    \bigg]^{1/2}. \label{eqn: y1 cauchy}
\end{align}
Noting that the the sum in (\ref{eqn: y1 cauchy}) is (almost) $\|Y^* Y\|_F^2$
and once again using the Cauchy--Schwarz inequality,
\begin{align}
   \sum_{i=1}^n \sum_{\substack{j=1 \\ j \neq i}}^n 
    |\langle y_i, y_j \rangle|^2 & 
    = \|Y^* Y\|_F^2 - \sum_{i=1}^n \|y_i\|_2^4 \nonumber \\
    ~ & = \mathrm{tr}(Y^* Y Y^* Y)
    - n \cdot \frac{1}{n} \sum_{i=1}^n (\|y_i\|_2^2)^2 \nonumber \\
    ~ & \leq \frac{n^2}{k} - n\left(\frac{1}{n}\sum_{i=1}^n 
    \|y_i\|_2^2 \right)^2 \\
    ~ & = \frac{n^2}{k} - n. \label{eqn: partial ineq}
\end{align}
Putting this back into (\ref{eqn: y1 cauchy}) we obtain the inequality
\begin{equation}
    \|Y^* Y\|_1 \leq n + \left[n(n-1)\left(\frac{n^2}{k} - n\right) 
    \right]^{1/2}.
\end{equation}
Equality is achieved in the Cauchy--Scwharz inequality
if and only if the vectors are scalar multiples.
For the first instance of Cauchy--Schwarz, this occurs if and only if $Y$ 
is equiangular. For the second instance of Cauchy--Schwarz, equality is
achieved if and only if $\|y_i\|_2^2$ is constant, that is $Y \in N(k,n)$.
Thus, equality is achieved in (\ref{eqn: welch y}) if and only if $Y$
is an equiangular tight frame. \qed
\end{prf}

Equality in (\ref{eqn: welch y}) depends on the existence of equiangular 
tight
frames of $n$ vectors in $\C^k$. The Gerzon bound says that if $Y$ forms an
equiangular tight frame of $n$ vectors in $\C^k$, then 
$n \leq k^2$~\cite{FicMix2015}. This gives the bound $k \geq 1/2 + 
\sqrt{1+4d}/2$ as a necessary condition for $Y$ to be an equiangular
tight frame. On the
other hand, if $Y$ in Lemma 1 is an equiangular tight frame, then $X$ is also an
equiangular tight frame since $X$ and $Y$ are Naimark complements
\cite{CasFickMixPet2013}. In particular, this gives the upper bound 
$k \leq d^2 - d$ as a necessary condition for $Y$ to be equiangular, 
because the Gerzon bound applied to $X$ gives the requirement that $n \leq d^2$.

Being within the range $1/2 + \sqrt{1+4d}/2 \leq k \leq d^2 - d$ is not 
a sufficient condition. The existence of equiangular tight frames of $n$
vectors in $k$ dimensions for $k+1 \leq n \leq k^2$ is an open question. 
Some equiangular tight
frames are known to exist, by construction, for certain values of $n$ and $k$.
An overview of the known constructions can be found in
\cite{FicMix2015}. On the other hand, there are known values of $n$ and $k$
for which equiangular tight frames cannot exist. One such example is the
case when $n = 8$ and $k = 3$ \cite{Szo2014}. In particular, equality in
the Welch bound is achieved for some values of $n$ and $k$ which satisfy 
$k+1 \leq n \leq k^2$, but not necessarily achieved at all values of $n$ 
and $k$ which satisfy that inequality.

By combining Lemma \ref{thm: minmax coherence} with Theorem \ref{lemma: welch y1}, we obtain

\begin{cor}[Welch]
Let $n > d$. For all $X \in N(d,n)$,
\begin{equation}
    \mu(X) \geq \sqrt{\frac{n-d}{d(n-1)}}.
\end{equation}
\end{cor}

\section{Bukh--Cox Bound via Linear Programming}
We now turn our attention to the range $1 \leq k < 1/2 + \sqrt{1+4d}/2$.
Since the Gerzon bound prevents $Y$ from being an equiangular tight
frame in this range, equality in (\ref{eqn: welch y}) cannot be achieved 
and a 
sharper bound can be obtained.
The Bukh--Cox bound is an improvement in this range, and is sharp if $Y$ is
given by 
concatenated copies of $k^2$ vectors in $\C^k$ which forms an equiangular 
tight frame.
In order to apply Delsarte's linear programming ideas, we require the following
special polynomials \cite{DelGoeSei1991}:
\begin{align*}
Q_0(x) & = 1, \\
Q_1(x) & = x - \frac{1}{k}, \\
Q_2(x) & = x^2 - \frac{4}{k+2} x + \frac{2}{(k+1)(k+2)}.
\end{align*}

\begin{thm}\thlabel{lemma: bukh-cox y1}
    For all $Y \in T(k,n)$ we have
\begin{equation}\label{eq: bukh-cox}
    \|Y^* Y\|_1 \leq \frac{n^2(1 + (k-1)\sqrt{1+k})}{k^2}.
\end{equation}
    Equality is
    achieved when $Y$ is of the form $Y = [Z | Z | \cdots | Z]$ where
    $Z \in \C^{k \times k^2}$ is an equiangular tight frame.
\end{thm}

\begin{figure*}
\centering
\includegraphics[width=0.95\textwidth]{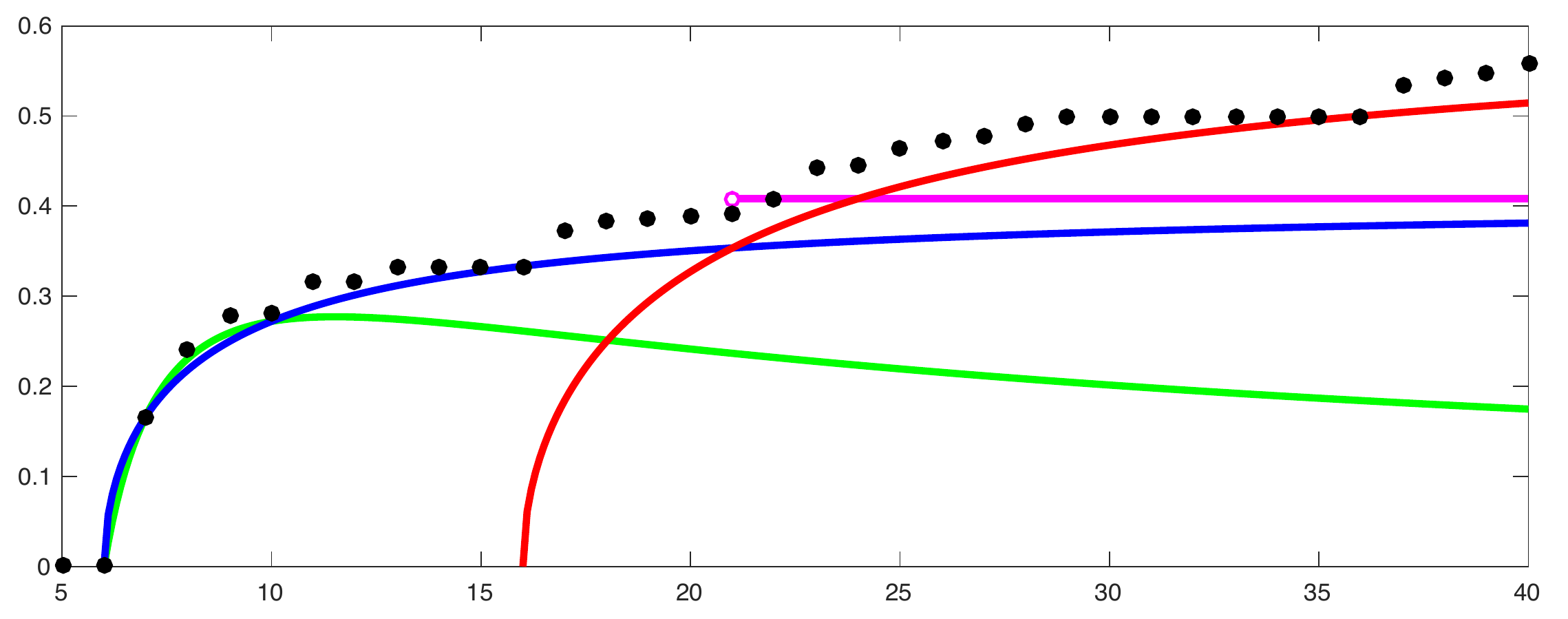}
\caption{The coherence of best known line packings in $\mathbb{R}^6$ for $n\in\{5,\ldots,40\}$, along with the best known lower bounds.
The black dots correspond to packings found in Sloane's database~\cite{sloane}. The Bukh--Cox bound is displayed in green, the Welch bound in blue, the orthoplex bound~\cite{conway1996packing} in pink, and the Levenshtein bound~\cite{haas2017levenstein} in red.
In this setting, the Bukh--Cox bound is the best known lower bound for $n=\{8,9\}$.}
\label{fig: bound}
\end{figure*}

\begin{prf}
By continuity, we may assume $y_i \neq 0$ for every $i \in \{1,
\cdots,n\}$ without loss of generality.
First, we normalize the columns of $Y$, $\{y_i\}_{i=1}^n$, by defining
$z_i := y_i/\|y_i\|_2$.
We obtain the desired bound considering the following linear program, 
inspired by Delsarte's LP bound,
\begin{align*}
\text{minimize } & c_0 \\
\text{subject to } & f(x) = c_0 Q_0(x) + c_1 Q_1(x) + c_2 Q_2(x) \\
~ & 0 \leq c_1 \leq k c_0, \,\, c_2 \leq 0, \\
~ & f(x) \geq \sqrt{x}, \,\, \forall x \in [0,1].
\end{align*}
Suppose we have a feasible $(c_0, c_1, c_2)$. Then,
\begin{align}
    \|Y^* Y\|_1 & = \sum_{i=1}^n \sum_{j=1}^n 
    |\langle z_i, z_j \rangle|
    \, \|y_i\|_2 \, \|y_j\|_2 \nonumber \\
    ~ & \leq \sum_{i=1}^n \sum_{j=1}^n
    f(|\langle z_i, z_j \rangle|^2) \|y_i\|_2 \, \|y_j\|_2 \nonumber \\
    ~ & = \sum_{\ell = 0}^2 \sum_{i=1}^n \sum_{j=1}^n 
    c_\ell Q_\ell(|\langle z_i, z_j \rangle|^2)\|y_i\|_2 \, \|y_j\|_2
    \nonumber \\
    ~ & = \sum_{\ell = 0}^2 c_\ell \sum_{i=1}^n \sum_{j=1}^n 
    Q_\ell(|\langle z_i, z_j \rangle|^2)\|y_i\|_2 \, \|y_j\|_2.
    \label{eqn: cl Ql}
\end{align}
We now establish an upper bound for each innermost term for 
$\ell \in \{0,1,2\}$.
For $\ell = 0$, since $Q_0(x) = 1$, we have
\begin{align}
~ & \sum_{i=1}^n \sum_{j=1}^n Q_0(|\langle z_i, z_j \rangle|^2) \|y_i\|_2 \,
\|y_j\|_2 &~ \nonumber \\
~ & = \sum_{i=1}^n \sum_{j=1}^n \|y_i\|_2 \, \|y_j\|_2 
= \left( \sum_{i=1}^n \|y_i\|_2 \right)^2
\nonumber \\
~ &=: S. \label{eqn: l0 bound}
\end{align}
For $\ell = 1$, we have
\begin{align}
~ & \sum_{i=1}^n \sum_{j=1}^n Q_1(|\langle z_i, z_j \rangle|^2) \|y_i\|_2 \,
\|y_j\|_2  \nonumber \\
~ & = \sum_{i=1}^n \sum_{j=1}^n \frac{|\langle y_i, y_j \rangle|^2}{\|y_i\|_2
\, \| y_j \|_2} - \frac{1}{k} \sum_{i=1}^n \sum_{j=1}^n \|y_i\|_2 \, \|y_i\|_2. \\
~ & = \sum_{i=1}^n \sum_{j=1}^n \frac{|\langle y_i, y_j \rangle|^2}{\|y_i\|_2
\, \| y_j \|_2} - S / k.
\label{eqn: l1 intermed}
\end{align}
To bound the first term of (\ref{eqn: l1 intermed}), 
we use Cauchy--Schwarz and the fact that $Y$ is an
$n/k$-tight frame:
\begin{align}
    ~ & \sum_{i=1}^n \sum_{j=1}^n \frac{|\langle y_i, y_j \rangle|^2}
    {\|y_i\|_2 \, \|y_j\|_2} \nonumber \\ 
    ~ & \leq \left( \sum_{i=1}^n \sum_{j=1}^n \frac{|\langle y_i, y_j
    \rangle|^2}{\|y_i\|_2^2} \right)^{\! 1/2}
    \left( \sum_{i=1}^n \sum_{j=1}^n \frac{|\langle y_i, y_j
    \rangle|^2}{\|y_j\|_2^2} \right)^{\! 1/2} \nonumber \\
    ~ & = \left( \sum_{i=1}^n \frac{n}{k} \frac{\|y_i\|_2^2}{\|y_i\|_2^2}
    \right)^{\! 1/2}
    \left( \sum_{j=1}^n \frac{n}{k} \frac{\|y_j\|_2^2}{\|y_j\|_2^2}
    \right)^{\! 1/2} = \frac{n^2}{k}.
\end{align}
Overall this gives the following bound for the $\ell = 1$ case:
\begin{equation}
\sum_{i=1}^n \sum_{j=1}^n Q_1(|\langle z_i, z_j \rangle|^2) \|y_i\|_2 \,
\|y_j\|_2 \leq \frac{1}{k}(n^2 - S) \label{eqn: l1 case}.
\end{equation}
Last, we need a bound for the $\ell = 2$ case. Let $\{e_m\}_{m=1}^{d_2}$ be
any orthonormal basis for the (finite) $d_2$-dimensional vector space spanned
by degree-$4$ projective harmonic polynomials 
in $k$ variables. Then, by the addition formula, there is a constant 
$C_{d_2,k}$, which depends on $d_2$ and $k$, such that
\begin{align}
~ & \sum_{i=1}^n \sum_{j=1}^n Q_2(|\langle z_i, z_j \rangle|^2) \|y_i\|_2 \,
\|y_j\|_2  \nonumber \\
~ & = C_{d_2,k} \sum_{i=1}^n \sum_{j=1}^n \sum_{m=1}^{d_2}
e_m(z_i) \overline{e_m(z_j)} \|y_i\|_2 \, \|y_j\|_2  \nonumber \\
~ & = C_{d_2,k} \sum_{m=1}^{d_2} \left| \sum_{i=1}^n e_m(z_i) \|y_i\|_2 
\right|^2 \geq 0. \label{eqn: l2 intermed}
\end{align}
Multiplying both sides of (\ref{eqn: l2 intermed}) by $c_2
\leq 0$ then gives
\begin{equation}\label{eqn: l2 final}
\sum_{i=1}^n \sum_{j=1}^n c_2Q_2(|\langle z_i, z_j \rangle|^2) \|y_i\|_2 \,
\|y_j\|_2 \leq 0.
\end{equation}
Finally, returning to (\ref{eqn: cl Ql}) we have
\begin{align}
\|Y^* Y\|_1 & = \sum_{\ell = 0}^2 c_\ell \sum_{i=1}^n \sum_{j=1}^n 
Q_\ell(|\langle z_i, z_j \rangle|^2)\|y_i\|_2 \, \|y_j\|_2 \nonumber \\
~ & \leq c_0 S + c_1\frac{1}{k}(n^2-S) 
= \left(c_0-\frac{c_1}{k}\right)S + c_1\cdot\frac{n^2}{k} \nonumber \\
~ & \leq c_0 n^2,
\label{eqn: c0 y bound}
\end{align}
where we have used that $S \leq n^2$. The bound~\eqref{eq: bukh-cox} comes from observing that the following choice of $(c_0, c_1, c_2
)$ is feasible:
\begin{align*}
c_0 & = \frac{1 + (k-1)\sqrt{1+k}}{k^2} \\
c_1 & = \frac{\sqrt{1+k} (-4 + k^2 + 4\sqrt{1+k})}{2k (2+k)} \\
c_2 & = \frac{-(2+4k+2k^2) + \sqrt{1+k}(2+3k+k^2)}{2k^2}.
\end{align*}
This feasible choice comes from forcing 
$f(\frac{1}{k+1}) = 1/\sqrt{k+1}$, $f(1) = 1$, and $f'(\frac{1}{k+1}) 
= \sqrt{k+1}/2$, and solving for $(c_0,c_1,c_2)$.
Equality is achieved in inequalities (\ref{eqn: cl Ql}), (\ref{eqn: l1 case}),
(\ref{eqn: l2 final}), (\ref{eqn: c0 y bound}) when
\begin{enumerate}
    \item $|\langle y_i, z_j\rangle| = |\langle z_i, y_j\rangle|,
    \,\, \forall i,j$,
    \item $f(|\langle z_i, z_j\rangle|^2) = |\langle z_i, z_j \rangle|, 
    \,\, \forall i,j$,
    \item $\sum_i \sum_j Q_2(|\langle z_i, z_j \rangle|^2) \|y_i\|_2 \,
    \|y_j\|_2 = 0$,
    \item $\|y_i\|_2 = 1, \,\, \forall i$.
\end{enumerate}
All four are achieved if $Y$ is multiple copies of an equiangular tight
frame of $k^2$ vectors in $\C^k$. \qed
\end{prf}

The proof of \thref{lemma: bukh-cox y1} actually
generates a bound for any feasible $(c_0, c_1, c_2)$ in the described
linear program. Minimizing $c_0$ gives the best possible bound
generated by this method. Computational experiments suggest that this particular
feasible $(c_0, c_1, c_2)$ gives the minimum $c_0$. 
Although equality
is achieved when $Y$ is multiple copies of an equiangular tight frame of $k^2$
vectors in $\C^k$, the existence of such frames is an open question, and
is known as Zauner's conjecture \cite{Zau1999}. 

By combining Lemma \ref{thm: minmax coherence} with Theorem \ref{lemma: bukh-cox y1}, we obtain

\begin{cor}[Bukh--Cox]
Let $n > d$. 
For all $X \in N(d,n)$,
\[
\mu(X) \geq \frac{(n-d)^2}{n+(n^2-nd-n)\sqrt{1+n-d} - (n-d)^2}.
\]
\end{cor}

Bukh and Cox also provide a new bound for the case of $n$ vectors in $\R^d$. For
the real case, it suffices to adjust the special polynomials in the proof of \thref{lemma: bukh-cox y1} \cite{DelGoeSei1991}. $Q_0(x)$ and $Q_1(x)$ remain the same,
but $Q_2(x)$ is replaced with:
\begin{equation*}
Q_2(x) = x^2 - \frac{6}{d+4} x + \frac{3}{(d+2)(d+4)}.
\end{equation*}
This adjustment changes the feasible region for the linear program and leads
to a different optimal $(c_0, c_1, c_2)$, and thus a different bound for the
$\R^d$ case. Fig.~\ref{fig: bound} demonstrates the Bukh--Cox bound
improvement over the Welch bound for small values of $k$ in the case where
the vectors are in $\R^6$. (We illustrate the real case since, in this case,
packing data is available and provided in \cite{sloane}.)

\section*{Acknowledgements}
We thank the anonymous reviewers for helpful comments that improved the 
presentation of our results.
MM and DGM were partially supported by AFOSR FA9550-18-1-0107. DGM was also supported by NSF DMS 1829955 and the Simons Institute of the Theory of Computing.

\bibliographystyle{IEEEtran}
\bibliography{sampta19}
\end{document}